\documentclass[12pt]{article}
\usepackage[a4paper]{geometry}
\usepackage{hyperref}
\usepackage{color,amsmath,amssymb,amsthm}
\usepackage{graphicx}

\usepackage{enumitem}

\newtheorem{theorem}{Theorem}[section]

\newtheorem{lemma}[theorem]{Lemma}
\newtheorem{definition}[theorem]{Definition}

\theoremstyle{remark}
\newtheorem{remark}[theorem]{Remark}

\definecolor{myurlcolor}{rgb}{0.6,0,0}
\hypersetup{colorlinks,
linkcolor=myurlcolor,
citecolor=myurlcolor,
urlcolor=myurlcolor}
\def\XXint#1#2#3{{\setbox0=\hbox{$#1{#2#3}{\int}$ }
\vcenter{\hbox{$#2#3$ }}\kern-.6\wd0}}

\theoremstyle{definition}

        \newcommand{\be}{\begin{equation}}
        \newcommand{\ee}{\end{equation}}
        \newcommand{\ba}{\begin{eqnarray}}
        \newcommand{\ea}{\end{eqnarray}}
        \newcommand{\ban}{\begin{eqnarray*}}
        \newcommand{\ean}{\end{eqnarray*}}
        \newcommand{\barr}{\begin{array}}
        \newcommand{\earr}{\end{array}}

\renewcommand{\to}{\rightarrow}

\newcommand{\maps}{\colon}

\newcommand{\R}{{\mathbb R}}

\date{}

\def\cH{\mathcal{H}}

\def\part{\partial}

\def\XXint#1#2#3{{\setbox0=\hbox{$#1{#2#3}{\int}$ }
\vcenter{\hbox{$#2#3$ }}\kern-.6\wd0}}

\title {Stability of Current Density Impedance Imaging}

\author{{Robert  Lopez \footnote{Department of Mathematics, University of California, Riverside, California, USA. E-mail: rlope021@ucr.edu.   }
\qquad Amir Moradifam\footnote{Department of Mathematics, University of California, Riverside, California, USA. E-mail: amirm@ucr.edu. Amir Moradifam is supported by NSF grant DMS-1715850.}}}

\begin{document}

\maketitle
{\small \noindent
%{\bf Keywords}  \\
%{\bf AMS subject classification} .
}
\begin{abstract}
We study stability of reconstruction in current density impedance imaging  (CDII), that is, the inverse problem of recovering the conductivity of a body from the measurement of the magnitude of the current density vector field in the interior of the object. Our results show that CDII is stable with respect to errors in interior measurements of the current density vector field, and confirm the stability of reconstruction which was previously observed in numerical simulations, and was long believed to be the case. 
\end{abstract}
\vskip 5em
\section{Introduction}
The classical Electrical Impedance Tomography (EIT) aims to obtain quantitative information on the electrical conductivity $\sigma$ of a conductive body  from measurements of voltages and corresponding currents at its boundary. Mathematics of EIT has been extensively studied, and many interesting results have been obtained about uniqueness, stability  and reconstruction algorithms for this problem. See  \cite{borcea, CI, isaacsonReview, GLKU}  for excellent reviews of the results. It is well known that 
that EIT is severely ill-posed, and provides images with very low resolution away from the boundary \cite{isacson, Man}.

A more recent class of Inverse Problems seeks to provide images with high accuracy and by using data obtained from the interior of the region. Such methods are referred to as Hybrid Inverse Problems or Coupled-physics methods, as they usually involve the interaction of two kinds of physical fields. In this paper we study stability of reconstruction in Current Density Impedance Imaging (CDII), that is, the inverse problem of recovering the conductivity of a body from the measurement of the magnitude of the current density vector field in the interior of the object. Interior measurements of current density is possible by Magnetic Resonance Imaging (MRI) due to the work of M. Joy  and his collaborators \cite{joy89, joy}. This problem has been  studied in \cite{MNT, MNTCalc, NTT07, NTT08, NTT10}. See also \cite{NTT11} for a comprehensive review. While the uniqueness of the reconstruction in CDII is established and a robust  reconstruction algorithm is developed in \cite{MNTim}, the stability of CDII is still open. In this paper, we aim to settle the stability of reconstruction in CDII, and provide a detailed stability analysis.

Let $\sigma$ be the isotropic conductivity of an object $\Omega\subset \R^n$,  $n\geq 2$, where $\Omega$ is a bounded open region in with connected boundary.  Suppose $J$ is the current density vector field generated by imposing a given boundary voltage $f$ on $\partial \Omega$.  Then the corresponding voltage potential $u$ satisfies the elliptic equation
\begin{equation}\label{PDE}
\nabla \cdot \left(\sigma \nabla u \right)=0, \ \ u|_{\partial \Omega}
=f.\end{equation}
By Ohm's law $ J=-\sigma \nabla u$, and $u$ is the unique minimizer of the weighted least gradient problem 
\begin{equation}\label{LGP}
I(w) =\min _{w\in BV_f(\Omega)} \int_{\Omega} a|\nabla w| dx,
\end{equation}
where $a=|J|$, and $BV_f(\Omega)= \{ w \in BV(\Omega), w|_{\partial \Omega} = f\}$, see \cite{MNT, MNTCalc, NTT07, NTT08, NTT10}.\\

Note that the weighted least gradient problem \eqref{LGP} is not strictly convex, and hence in general it may not have a unique minimizer. See \cite{JMN} where the second author and his collaborators showed that for $a\in C^{1,\alpha} (\Omega)$, $0<\alpha <1$, the least gradient problem \eqref{LGP} could have infinitely many minimizers. Since any stability result trivially implies uniqueness, it is evident that one needs additional assumptions to prove any stability result. Indeed stability analysis of CDII is a challenging problem. The first stability result on CDII was proved by Montalto and Stefanov in \cite{MS}. 

\begin{theorem}[\cite{MS}]
Let $u$ solve equation (1) and let $\tilde{u}$ solve equation (1) for $\tilde{\sigma}$ with $|\nabla \tilde{u}| > 0$ in $\overline{\Omega}$. For any $0 < \alpha < 1$, there exists $s > 0$ such that if $\| \sigma \|_{H^s (\Omega)} < L$ for some $L > 0$ then there is an $\epsilon > 0$ such that if 
\begin{equation}\label{epsilon}
 \| \sigma - \tilde{\sigma} \|_{C^2(\Omega)} < \epsilon,
\end{equation}
then 
\[ \| \sigma - \tilde{\sigma} \|_{L^2(\Omega)} < C \||J| - |\tilde{J}| \|_{L^2(\Omega)}^{\alpha} \]

\end{theorem}

\noindent Later in  \cite{MT}, Montalto and Tamasn proved the following stability result.

\begin{theorem}[\cite{MT}]
Let $\sigma \in C^{1,\alpha} (\overline{\Omega})$, $0 < \alpha < 1$, be positive in $\overline{\Omega}$. Let $u$ solve equation (1) with $|\nabla u| > 0$ in $\overline{\Omega}$. There exists $\epsilon > 0$ depending on $\Omega$ and some $C > 0$ depending on $\epsilon$ such that if $\tilde{\sigma} \in C^{1,\alpha} (\overline{\Omega})$ with $\tilde{u}$ solving (1) for $\tilde{\sigma}$, $u = \tilde{u} = f$ on $\partial \Omega$,  $\sigma = \tilde{\sigma}$ on $ \partial \Omega$, and 
\[ \| \sigma - \tilde{\sigma} \|_{C^{1,\alpha} (\overline{\Omega})} < \epsilon,  \]
then 
\[ \| \sigma - \tilde{\sigma} \|_{L^2 (\Omega)} \leq C \| \nabla \cdot (\Pi_{\nabla u} (J - \tilde{J})) \|_{L^2 (\Omega)}^{\frac{\alpha}{2 + \alpha}} \]
where $\Pi_{\nabla u} (J - \tilde{J})$ is the projection of $J - \tilde{J}$ onto $\nabla {u}$.
\end{theorem}

Note that both of the above results assume a priori that $\sigma$ and $\tilde{\sigma}$ are close, and a natural question which remains open is that whether there exists two distant conductivities $\sigma$ and $\tilde{\sigma}$ which could induce the corresponding currents $J$ and $\tilde{J}$  with $| |J|-|\tilde{J}| |$ arbitrarily small.  In this paper we address the this question and show that the answer is negative, and hence show that CDII is actually stable. Under some natural assumption, we shall prove that in dimensions $n=2,3$ the following stability result holds

\begin{equation}
 \| \sigma-\tilde{\sigma}  \|_{L^1 (\Omega)} \leq C\||J| - |\tilde{J}|\|^{\frac{1}{4}}_{L^{\infty}({\Omega})},
\end{equation}
for some constant $C$ independent of $\tilde{\sigma}$ (see Theorems \ref{sigmaStability2} and \ref{sigmaStability3} for precise statements of the results). 

The paper is organized as follows. In Section 2, under very weak assumptions, we will prove that the structure of level sets of the least gradient problem \eqref{LGP} is stable. In Section 3, we will provide stability results for minimizers of \eqref{LGP} in $L^1$. In Section 4, we will prove stability of minimizers of \eqref{LGP} in $W^{1,1}$, and shall use them to prove Theorems \ref{sigmaStability2} and \ref{sigmaStability3} which are the main results of this paper.

\section{Stability of level sets}
In this section, we show that the structure of the level sets of minimizers of the least gradient problem \eqref{LGP} is stable. Throughout the paper, we will assume that $a, \tilde{a} \in C (\Omega)$ with 
\begin{equation} \label{mM}
0<m \leq  a(x),\tilde{a}(x) \leq M, \ \ \forall x\in \Omega,
\end{equation}
for some positive constants $m,M$. The following theorem which was proved in \cite{Struc1} by the second author, shall play a crucial role in the proof of the results in this section. 

\begin{remark} \label{remark1}
In general the least gradient problem \eqref{LGP} may not have a minimizer \cite{Dos, ST}. Throughout the paper we shall assume that \eqref{LGP} has a solution. For sufficient conditions for the existence of minimizers of weighted least gradient problems we refer to \cite{G, JMN, Struc1}. Note also that any voltage potential $u$ solving the equation \eqref{PDE} is also a minimizer of \eqref{LGP}. In particular, if $0<a(x)\in C(\overline{\Omega})$ and $\partial \Omega$ satisfies a Barrier condition (see Definition 3.1 in \cite{JMN}), then for every $f \in C(\partial \Omega)$ the least gradient problem \eqref{LGP} has a minimizer in $BV_f(\Omega)$.  In other words the set of weights for which the least gradient problem \eqref{LGP} has a solution is open in $C(\overline{\Omega})$ if  $\partial \Omega$ satisfies a barrier condition. 
\end{remark}

\begin{theorem}[\cite{Struc1}]\label{Structure1}
Let $\Omega \subset \R^n$ be a bounded open set with Lipschitz boundary and assume that $a \in C(\overline{\Omega}) $ is a non-negative function, and $f\in L^1(\partial \Omega)$.  Then there exists a divergence free vector field $J \in (L^{\infty}(\Omega))^n$ with $|J|\leq a$ a.e. in $\Omega$ such that  every minimizer $w$ of 
\eqref{LGP}
 satisfies 
\begin{equation}\label{parallel}
J\cdot \frac{Dw}{|Dw|}=|J|=a, \ \ |Dw|-a.e. \ \ \hbox{in} \ \ \Omega,
\end{equation}
where $\frac{Dw}{|Dw|}$ is the Radon-Nikodym derivative of $Dw$ with respect to $|Dw|$. 

\end{theorem}

\begin{lemma} \label{lemma1} Let $f\in L^{1}(\partial \Omega)$, and assume $u$ and $\tilde{u}$ are minimizers of \eqref{LGP} with the weights $a$ and $\tilde{a}$, respectively.  Then
\begin{equation}\label{first}
  \left| \int_{\Omega} a|D u| dx - \int_{\Omega} \tilde{a}|D \tilde{u}| dx \right| \leq C \| a - \tilde{a} \|_{L^{\infty}(\Omega)},
\end{equation}
for some constant $C=C(m,M,\Omega, f)$ independent of $u$ and $\tilde{u}$.
\end{lemma}

\vskip 1em \noindent
\textbf{Proof. } First note that in view of \eqref{mM} we have
\[ m \int_{\Omega} |D \tilde{u}| dx \leq \int_{\Omega} \tilde{a} |D \tilde{u}| dx \leq \int_{\Omega} \tilde{a} |D w| dx \leq M \int_{\Omega}|D w|  \]
for any $w \in BV_f(\Omega)$. Thus $ \int_{\Omega}|D\tilde{u}|\leq C$, and similarly $\int_{\Omega}|Du|\leq C$ for some constant $C$ which depends only on $m,M,$ and $\Omega$. Hence
\begin{equation}\label{C}
\max \left\{ \int_{\Omega}|D\tilde{u}|, \int_{\Omega}|Du| \right\}\leq C,
\end{equation}
for some $C(m,M)$ independent of $\tilde{u}$ and $u$. Since $u$, $\tilde{u}$ are the  minimizers of \eqref{LGP} with the weights $a$ and $\tilde{a}$, 

\begin{align*}
 \int_{\Omega} a|D u| dx - \int_{\Omega} \tilde{a}|D u| dx  & \leq \int_{\Omega} a|D u| dx - \int_{\Omega} \tilde{a}|D \tilde{u}| dx\\
 & \leq   \int_{\Omega} a|D \tilde{u}| dx - \int_{\Omega} \tilde{a}|D \tilde{u}| dx.
\end{align*}
Thus 
\[   \int_{\Omega} (a - \tilde{a}) | D u | dx \leq \int_{\Omega} a|D u| dx - \int_{\Omega} \tilde{a}|D u| dx \leq \int_{\Omega} (a - \tilde{a}) | D \tilde{u} | dx, \]
and we get
\begin{align*}
   -\| a - \tilde{a} \|_{L^{\infty}(\Omega)} \| D u \|_{L^{1}(\Omega)} & \leq \int_{\Omega} a|D u| dx - \int_{\Omega} \tilde{a}|D u| dx \\ 
   & \leq \| a - \tilde{a} \|_{L^{\infty}(\Omega)} \| D \tilde{u} \|_{L^{1}(\Omega)} .
\end{align*}
Hence \eqref{first} follows from \eqref{C}. $\hfill$ $\square$
\vspace{.5cm}

Let $\nu_\Omega$ denote the outer unit normal vector to $\partial\Omega$. Then for every $T\in(L^\infty(\Omega))^n$ with $\nabla\cdot T\in L^n(\Omega)$ there exists a unique function $[T,\nu_\Omega]\in L^\infty(\partial\Omega)$ such that 
\begin{equation}\label{trace}
\int_{\partial\Omega}[T,\nu_\Omega]u\,d\cH^{n-1}=\int_\Omega u\nabla\cdot Tdx+\int_\Omega T\cdot D udx,
\quad u\in C^1(\bar\Omega).
\end{equation}
Moreover, for $u\in BV(\Omega)$ and $T\in( L^\infty(\Omega))^n$ with $\nabla \cdot T\in L^n(\Omega)$, the linear functional $u\mapsto(T\cdot Du)$ gives rise to a Radon measure on $\Omega$, and \eqref{trace} holds for every $u\in BV(\Omega)$ (see \cite{Alberti, Anzellotti} for a proof). We shall need the weak integration by parts formula \eqref{trace}.

\vskip 1em \noindent
\begin{lemma}\label{lemma2}
 Let $f\in L^1(\partial \Omega)$, and assume $u$ and $\tilde{u}$ are minimizers of \eqref{LGP} with the weights $a$ and $\tilde{a}$, respectively. Let $J$ and $\tilde{J}$ be the divergence free vector fields guaranteed by Theorem \ref{Structure1}. Suppose $0 \leq \sigma (x) \leq \sigma_1$ in $\Omega$ for some constant  $\sigma_1>0$, where $\sigma$ is the Radon-Nikodym derivative of $|J|dx$ with respect to $|Du|$ . Then 
\begin{equation}
 \int_{\Omega} |J||\tilde{J}| - J \cdot \tilde{J} dx   \leq C \| a - \tilde{a} \|_{L^{\infty}(\Omega)}, 
\end{equation}
where $C=C(m,M,\sigma_1, \Omega, f, u)$ is a constant independent of $\tilde{a}$. 
\end{lemma}
\vskip 1em \noindent
\textbf{Proof.} We have 
\begin{align*}
  \int_{\Omega} |J||\tilde{J}| - J \cdot \tilde{J} dx & =  \int_{\Omega}  \sigma |\tilde{J}| |D u|  - \sigma \tilde{J} \cdot D u  dx \\
  & \leq \sigma_1 \int_{\Omega}  |\tilde{J}| |D u|  - \tilde{J} \cdot D u  dx \\
  & = \sigma_1 \int_{\Omega}  |\tilde{J}| |D u| dx  -\int_{\partial \Omega} f [\tilde{J},\nu_{\Omega}]  dx\\
  & = \sigma_1 \int_{\Omega}  |\tilde{J}| |D u|  -\tilde{J} \cdot D \tilde{u}  dx\\
  & = \sigma_1 \int_{\Omega}  |\tilde{J}| |D u|  - |\tilde{J}| | D \tilde{u}|  dx,
\end{align*}
where we have used \eqref{parallel} and the integration by parts formula \eqref{trace}. On the other hand it follows from lemma \ref{lemma1} that 
\begin{align*}
  \sigma_1 \int_{\Omega}  |\tilde{J}| |D u|  - |\tilde{J}| |D \tilde{u}|  dx & =  \sigma_1 \int_{\Omega}  |\tilde{J}| |D u|  - |J||D u| + |J||D u| - |\tilde{J}| |D \tilde{u}|  dx  \\
   & = \sigma_1 \left( \int_{\Omega}  (a - \tilde{a}) |D u| dx + \int_{\Omega} a |D u|  - \tilde{a} |D \tilde{u}| dx \right)   \\
   & \leq  \sigma_1 (\|D u\|_{L^1(\Omega)}  \| a - \tilde{a} \|_{L^{\infty}(\Omega)} + C \| a - \tilde{a} \|_{L^{\infty}(\Omega)} ),
\end{align*}
which yields the desired result. \hfill  $\square$

\vskip 1em \noindent
Roughly speaking, Lemma \ref{lemma2} implies that as $a \to \tilde{a}$,  $\frac{D u}{|Du|}(x)$ becomes parallel to $\frac{D\tilde{u}}{|D\tilde{u}|} (x)$ at points where the two gradients do not vanish. We are now ready to prove the main result of this section. \\

\begin{theorem}\label{levelSetTheorem}
Let $f\in L^1(\partial \Omega)$, and assume $u$ and $\tilde{u}$ are minimizers of \eqref{LGP} with the weights $a$ and $\tilde{a}$, respectively. Let $J$ and $\tilde{J}$ be the divergence free vector fields guaranteed by Theorem \ref{Structure1}. Suppose $0 \leq \sigma (x) \leq \sigma_1$ in $\Omega$ for some constant  $\sigma_1>0$, where $\sigma$ is the Radon-Nikodym derivative of $|J|dx$ with respect to $|Du|$ . Then 
\begin{equation}\label{JsAreClose}
 \|J - \tilde{J} \|_{L^1(\Omega)} \leq C \|a - \tilde{a}\|^{\frac{1}{2}}_{L^{\infty}(\Omega)},
\end{equation}
where $C=C(m,M,\sigma_1, \Omega, f, u)$ is a constant independent of $\tilde{a}$. 
\end{theorem}

\vskip 1em \noindent
\textbf{Proof.} We have 
\begin{align*}
 \sqrt{|J - \tilde{J}|^2} & = \sqrt{|J|^2 + |\tilde{J}|^2 - 2J\cdot \tilde{J}} \\
     & = \sqrt{\left||J| - |\tilde{J}| \right|^2 + 2(|J||\tilde{J}| - J\cdot \tilde{J})} \\
     & \leq \left||J| - |\tilde{J}| \right| + \sqrt{2(|J||\tilde{J}| - J\cdot \tilde{J})}.
\end{align*}
Hence,
\begin{align*}
 \|J - \tilde{J} \|_{L^1(\Omega)} & = \int_{\Omega} \sqrt{\left| |J| - |\tilde{J}| \right|^2} dx \\
   & \leq \int_{\Omega} \left||J| - |\tilde{J}| \right| dx + \int_{\Omega} \sqrt{2(|J||\tilde{J}| - J\cdot \tilde{J})} dx  \\
    & = \int_{\Omega} |a - \tilde{a}| dx + \int_{\Omega} \sqrt{2(|J||\tilde{J}| - J\cdot \tilde{J})} dx\\
      &\leq |\Omega| \|a - \tilde{a}\|_{L^{\infty}(\Omega)} + |\Omega|^{1/2} \left(\int_{\Omega} 2(|J||\tilde{J}| - J\cdot \tilde{J}) dx \right)^{1/2}   \\
      & \leq |\Omega| \|a - \tilde{a}\|_{L^{\infty}(\Omega)} + (2  |\Omega|)^{1/2} ( C \| a - \tilde{a} \|_{L^{\infty}(\Omega)} )^{1/2} \\
      &= \left(|\Omega| \|a - \tilde{a}\|^{\frac{1}{2}}_{L^{\infty}(\Omega)} +(2 |\Omega|)^{1/2} C^{\frac{1}{2}}\right) \|a - \tilde{a}\|^{\frac{1}{2}}_{L^{\infty}(\Omega)},
\end{align*}
where we have used the Holder's inequality and Lemma \ref{lemma2}. \hfill  $\square$
\vskip 1em 
\begin{remark} In view of Theorem \ref{Structure1}, $\frac{D u}{|Du|}$ and $\frac{D \tilde{u}}{|D\tilde{u}|}$ are parallel to $J$ and $\tilde{J}$, respectively.  So Theorem \ref{levelSetTheorem} implies that  if $\tilde{a}$ is close to $a$, then the structure of level sets of $\tilde{u}$ is close to that of $u$. 
\end{remark}

\section{$L^{1}$ stability of the minimizers}
In this section, we establish stability of minimizers of the least gradient problem \eqref{LGP} in $L^1$. In general \eqref{LGP} does not even have unique minimizers, so in order to prove any stability results further assumptions on the weights $a,\tilde{a}$, and on the corresponding minimizers are expected.  

\begin{definition} Fix the positive constants $\sigma_0, \sigma_1 \in \R$. We say that $u\in C^1(\bar{\Omega})$ is  admissible if it solves the conductivity equation \eqref{PDE} for some $\sigma \in C(\Omega)$ with 
\[0<\sigma_0<\sigma \leq \sigma_1,\]
and $m \leq |J|=|\sigma \nabla u| \leq M$, where $m$ and $M$ are positive constants as in \eqref{mM}. We shall denote the corresponding induced current by $J=-\sigma \nabla u$. 
\end{definition}

\begin{remark} \label{remark2}
Let  $\Omega \subset \R^n$ with $n \geq 2$ be a bounded Lipschitz domain and suppose $\partial \Omega$ satisfies the barrier condition defined in Definition 3.1 in \cite{JMN}). A. Zuniga proved in \cite{Zuniga} that if $0<a \in C^2({\overline{\Omega}})$, then for any boundary data $f \in C(\partial \Omega)$ the least gradient problem \eqref{LGP}
has a minimizer $u \in C(\overline{\Omega})$. If $|\nabla u|>0$ in $\overline{\Omega}$, then 
\[\sigma =\frac{a}{|\nabla u|} \in C(\overline{\Omega}),\]
and by elliptic regularity $u\in C^1(\overline{\Omega})$, and therefore \eqref{LGP} has an admissible minimizer. To guarantee the condition $|\nabla u|>0$ on $\overline{\Omega}$, in dimension $n=2$ it suffices to assume that the boundary data $f\in \partial \Omega$ is two-to-one, i.e. $f$ has only two critical points on $\partial \Omega$ (see Theorem 1.1 in \cite{Ales}). In higher dimensions, it is still an open problem to find sufficient conditions under which $|\nabla u|>0$ on $\overline{\Omega}$. 
\end{remark}

We will first prove our results in dimension $n=2$ and then extend them to dimensions $n=3$. 

Let $u\in C^1(\Omega)$ with $|\nabla u|>0$ in $\Omega$. Then it follows from the regularity result of De Giorgi (see, e.g, Theorem 4.11 in \cite{G}) that all level sets of $u$ are $C^{1}$ curves. We will assume that the length of level sets of $u$ in $\Omega$ is uniformly bounded, i.e. 
\begin{equation}\label{bpundedLength}
\sup_{t\in \R} \int_{\{u=t\}\cap \Omega}1 dl =L_M<\infty.\\
\end{equation}
\vspace{.2cm}

\begin{theorem}\label{usAreClose}
 Let $n=2$, and suppose $u$ and $\tilde{u}$ are admissible with $u|_{\partial \Omega}=\tilde{u}|_{\partial \Omega}=f,$ and corresponding current density vector fields $J$ and $\tilde{J}$, respectively. If $u$ satisfies \eqref{bpundedLength}, then 
 \begin{equation}
\|u - \tilde{u} \|_{L^1 (\Omega)} \leq C\parallel |J| - |\tilde{J}|\parallel^{\frac{1}{2}}_{L^{\infty}({\Omega})},
\end{equation}
for some constant $C(m,M,\sigma_0,\sigma_1,f, u, L_M)$ independent of $\tilde{u}$ and $\tilde{\sigma}$. 
\end{theorem}
\vskip 1em \noindent
\textbf{Proof.} Since $u$ is admissible, 
\[|\nabla u(x)|=\frac{|J(x)|}{\sigma(x)}\geq \frac{m}{\sigma_1}>0, \ \ \forall x\in \Omega.\]
Using the coarea formula we get 
\begin{equation}\label{leveleSetBound}
\frac{m}{\sigma_1} \int_{\Omega} |u - \tilde{u}| dx \leq \int_{\Omega} |\nabla u| |u - \tilde{u}| dx = \int_{\R} \int_{\{u = t\}\cap \Omega} |u - \tilde{u}|  dl dt.
\end{equation}
Since $|\nabla u|>0$ in $\Omega$, it follows from the regularity result of De Giorgi (Theorem 4.11 in \cite{G}) that all level sets of $u$ are $C^{1}$ curves. Now let $\Gamma_{t}$  be a connected component of $\{ x \in \Omega \maps u(x) = t \} \subset \Omega$, and  $\gamma \maps [0,L] \to \Gamma_{t}$ to be a path parameterized by the arc length with $\gamma(0)\in \partial \Omega$. Define 
\[h(s) := u(\gamma(s)) - \tilde{u} (\gamma(s)).\]
Then $h(0) = 0$. Moreover since  $\nabla u (\gamma(s)) \cdot \gamma'(s)=0$ on $\Gamma_t$,
\begin{eqnarray*}
h'(s)& = &\nabla u (\gamma(s)) \cdot \gamma'(s) - \nabla \tilde{u} (\gamma(s)) \cdot \gamma'(s)  \\
&=&\left( \frac{\sigma}{\tilde{\sigma}}(\gamma(s))\nabla u (\gamma(s) - \nabla \tilde{u} (\gamma(s) \right) \cdot \gamma'(s). 
\end{eqnarray*}
We can rewrite the above equality as
\[  h'(s) = \frac{J(\gamma(s))-\tilde{J}(\gamma(s))}{\tilde{\sigma}(\gamma(s))} \cdot \gamma'(s).\]
Now let $x^*_t$ be a point on $\Gamma_t$ where the maximum distance between $u$ and $\tilde{u}$ along the path $\gamma$ occurs, i.e. 
\[|u(x^*_t) -  \tilde{u}(x^*_t)|=\max \limits_{x \in \Gamma_t} |u(x) -  \tilde{u}(x)|.\]
Then $x^*_t=\gamma(s_0) $ for some $s_0\in [0,L]$, and 
\begin{eqnarray*}
|u(x^*_t) -  \tilde{u}(x^*_t)| = |h(s_0)| &=& \left| \int_{0}^{s_0} \frac{J(\gamma(\tau))-\tilde{J}(\gamma(\tau))}{\tilde{\sigma}(\gamma(\tau))} \cdot \gamma'(\tau) d\tau \right| \\
& \leq & \int_{0}^{s_0} \frac{1}{\tilde{\sigma}(\gamma(\tau))}|J(\gamma(\tau))-\tilde{J}(\gamma(\tau))| d\tau\\
&\leq &\frac{1}{\sigma_0} \int_{0}^{s_0} |J(\gamma(\tau))-\tilde{J}(\gamma(\tau))| d\tau.
\end{eqnarray*}
In particular for every $x \in \Gamma_t$
\[ |u(x)-\tilde{u}(x)| \leq |u(x^*_t) -  \tilde{u}(x^*_t)| \leq \frac{1}{\sigma_0}\int_{0}^{L} |J(\gamma(\tau))-\tilde{J}(\gamma(\tau))| d\tau, \]
where $L$ denotes the entire length of $\Gamma_t$.  Hence
\begin{eqnarray*}
\int_{\Gamma_t}|u(x)-\tilde{u}(x)|dl & \leq & |u(x^*_t) -  \tilde{u}(x^*_t)| \int_{\Gamma_t} 1 dl\\
& \leq & L_M  |u(x^*_t) -  \tilde{u}(x^*_t)| \\
& \leq & \frac{L_M}{\sigma_0}\int_{0}^{L} |J(\gamma(\tau))-\tilde{J}(\gamma(\tau))| d\tau \\
& =& \frac{L_M}{\sigma_0} \int_{\Gamma_t}|J-\tilde{J}|dl,
\end{eqnarray*}
and therefore
\begin{equation} \label{es1}
\int_{\{u = t\}\cap \Omega} |u - \tilde{u}|  dl \leq \frac{L_M}{\sigma_0} \int_{\{u = t\}\cap \Omega} |J-\tilde{J}| dl.
\end{equation}
Thus we have 

\begin{align*}
   \int_{\R} \int_{\{u = t\}\cap \Omega} |u - \tilde{u}|  dl dt &  \leq \frac{L_M}{\sigma_0} \int_{\R} \int_{\{u = t\}} |J - \tilde{J}| dl dt \\
      & = \frac{L_M}{\sigma_0} \int_{\Omega} |\nabla u| |J - \tilde{J}| dx \\
      & \leq \frac{L_M }{\sigma_0 }\|\nabla u\|_{L^{\infty}({\Omega})} \int_{\Omega} |J - \tilde{J}| dx  \\
     & \leq C\|a - \tilde{a}\|^{\frac{1}{2}}_{L^{\infty}({\Omega})}
\end{align*}
$C(m,M,\sigma_0,\sigma_1,f, u, L_M)$ independent of $\tilde{u}$ and $\tilde{\sigma}$, where we have used \eqref{es1} and Theorem \ref{levelSetTheorem}. \hfill $\square$
\vspace{.3cm}

Next we generalize Theorem \ref{usAreClose} to dimension $n=3$. In order to do this, we need the following additional assumption on level sets of $u$. 
\begin{definition} \label{foliation} Let $u\in C^1(\bar{\Omega})$ be admissible. We say that level sets of $u$ can be foliated to one-dimensional curves if for almost every $t \in range (u)$, every conected component $\Gamma_t$ of $\{u=t\}$ (equipped with the metric induced from the Euclidean metric in $\R^3$) there exists a function $g_t(x) \in C^1(\Gamma_t)$ such that $0<c_g \leq |\nabla g_t| \leq C_g$, for some constants $c_g$ and $C_g$ independent of $t$. Moreover, every connected component  of $\{u=t\}\cap \{g_t=r\}\cap \Omega$ is a $C^1$ curve reaching the boundary $\partial \Omega$ for almost every $t\in range(u)$ and all $r \in \R$. Similar to the case $n=2$, we assume that the length of connected components of  $\{u=t\}\cap \{g_t=r\}\cap \Omega$ are uniformly bounded by some constant $L_M.$
\end{definition}

\begin{remark}\label{foliationRemark} It follows from the regularity result of De Giorgi (see, e.g. Theorem 4.11 in \cite{G}) that for a function $u\in BV(\Omega)$, level sets $\{u=t\}$ is a $C^1$-hypersurface
for almost all $t\in range(u)$. Note also that every connected component of $\{u=t\}$ reaches the boundary $\partial \Omega$ (see \cite{MNT, MNTCalc, NTT07, NTT08}), for almost every $t$. Now let $\Gamma_t$ be a $C^1$ connected component of $\{u=t\}$. If $f$ has only two critical points (one minimum and one maximum points) on $\partial \Omega$, then $\Gamma_t$ is a simply-connected $C^1$ surface reaching the boundary $\partial \Omega$, and hence there exists a $C^1$ homeomorphism $\mathcal{F}_t$ from $\overline{B(0,1)} \subset \mathbb{R}^2$ to the closure of $\Gamma_t$ in $\overline{\Omega}$ (see Theorem 3.7 and Theorem 2.9 in \cite{Morris}). It is easy to see that the unit ball $B(0,1)$ can be foliated to one dimensional curves by level sets of $g:B(0,1) \rightarrow \mathbb{R}$ defined by $g(x,y)=y$. Consequently $\Gamma_t$ can be foliated into one dimensional curves reaching the boundary of $\partial \Omega$ by level sets of $g_t(X)=g(\mathcal{F}_{t}^{-1}(X))$, $X \in \Gamma_t$. Note also that since $g$ and $\mathcal{F}^{-1}_{t}$ are both $C^1$, and since $\overline{\Gamma_t}$ is compact, there exists constant $c(t), C(t)>0$ such that 
\begin{equation}\label{g_tBound}
0<c(t)<|\nabla g_t|<C(t)    \ \ \ \  \hbox{on}  \ \ \Gamma_t.
\end{equation}
Indeed the above argument shows that \eqref{g_tBound} holds for every connected components of almost every level sets of a function $u\in BV(\Omega)$, for some constant $c(t), C(t)$ depending on $t$. So in Definition \ref{foliation} the only significant assumption is that the constants $c(t)$ and $C(t)$ are uniformly bounded from below and above by two positive constant $c_g$ and $C_g$. In particular, if $u$ is a $C^1$ function with $|\nabla u|>0$ in $\Omega$ and $\{x\in \partial \Omega: f(x)=t\}$ has finitely many connected components for all $t$, then it follows from the implicit function theorem that every level set of $u$ is a $C^1$ surface, and hence existence of $c_g$ and $C_g$ follows immediately from compactness of $range(u)$, and hence level sets of $u$ can be foliated to one-dimensional curves in the sense of Definition \ref{foliation}.  

\end{remark}

\vskip 1em \noindent

\begin{definition}
Let $t\in range(u)$ and suppose $\Gamma_t^{i}$, $i\in I$, are $C^1$ connected components of $\{u=t\}$, where $I$ is countable. In view of Remark \ref{foliationRemark}, there exists functions $g_t^i:\Gamma_t^{i} \rightarrow R$ whose level sets foliate $\Gamma_t^{i}$ into one dimensional curves in the sense of Definition \ref{foliation}. We define $g_t : \{u=t\} \rightarrow R$ be the function with
\begin{equation}
g_t|_{\Gamma_t^i}=g_t^{i},  \ \ \ \ i\in I.
\end{equation}
We shall use this notation throughout the paper. 

\end{definition}

\begin{theorem}
Let $n=3$, and suppose $u$ and $\tilde{u}$ are admissible with $u|_{\partial \Omega}=\tilde{u}|_{\partial \Omega}=f$ and corresponding current density vector fields $J$ and $\tilde{J}$, respectively. Suppose the level sets of $u$ can be foliated to one-dimensional curves in the sense of Definition \ref{foliation}. Then
\begin{equation}
\|u - \tilde{u} \|_{L^1 (\Omega)} \leq C \||J| - |\tilde{J}| \|^{\frac{1}{2}}_{L^{\infty} (\Omega)},
\end{equation}
where $C(m,M,\sigma_0,\sigma_1,f, u, L_M, c_g,C_g)$ is independent of $\tilde{u}$ and $\tilde{\sigma}$. 
\end{theorem}
\vskip 1em \noindent
\textbf{Proof.} The proof is similar to the proof of Theorem \ref{usAreClose}, and we provide the details for the sake of the reader. Since $u$ is admissible, 
\begin{equation}\label{firstSt}
\frac{m}{\sigma_1} \int_{\Omega} |u - \tilde{u}| dx \leq \int_{\Omega} |\nabla u|  |u - \tilde{u}| dx = \int_{\R} \int_{\{u = t\}\cap \Omega} |u - \tilde{u}| dS dt. 
\end{equation}

The level sets of $u$ can be foliated into one-dimensional curves by level sets of some function $g_t$ in the sense of Definition \ref{foliation}. Thus 
\begin{align*}  \int_{\R} \int_{\{u = t\}\cap \Omega} |u - \tilde{u}| dS dt & =  \int_{\R} \int_{\{u = t\}\cap \Omega} \frac{|\nabla g_t|}{|\nabla g_t|} |u - \tilde{u}| dS dt \\
               & = \int_{\R} \int_{\R} \int_{\{u = t\} \cap \{g_t = r\}\cap \Omega} \frac{1}{|\nabla g_t|} |u - \tilde{u}| dl dr dt  \\
               & \leq \frac{1}{c_g} \int_{\R} \int_{\R} \int_{\{u = t\} \cap \{g_t = r\}\cap \Omega}  |u - \tilde{u}| dl dr dt.
 \end{align*}
Similar to the two dimensional case, we parameterize every connected component $\Gamma_t$ of $\{u = t\} \cap \{g_t = r\}\cap \Omega$ by arc length, $\gamma \maps [0, L] \to \Gamma_t$ with $\gamma(0) \in \partial \Omega$, and let $h(s) = u(\gamma(s)) - \tilde{u}(\gamma(s))$. Let $x^*_t$ be the point that maximizes $|u - \tilde{u}|$ on $\Gamma_t$ and suppose $\gamma(s_0) = x^*_t$  for some $s_0 \in (0,L)$, where $L$ is the length of $\Gamma_t$. Then by an argument similar to the one in the proof of Theorem \ref{usAreClose} we get 
\[   |u(x^*_t) - \tilde{u}(x^*_t)| \leq \frac{1}{\sigma_0}\int_{0}^{L} |J(\gamma(\tau)) - \tilde{J}(\gamma(\tau))| d\tau,\]
and consequently 

\begin{eqnarray*}
\int_{\Gamma_t}|u(x)-\tilde{u}(x)|dl & \leq & \frac{L_M}{\sigma_0} \int_{\Gamma_t}|J-\tilde{J}|dl.
\end{eqnarray*}
Hence
\begin{equation}
\int_{\{u = t\} \cap \{g_t = r\}\cap \Omega}  |u - \tilde{u}| dl \leq \frac{L_M}{\sigma_0} \int_{\{u = t\} \cap \{g_t = r\} \cap \Omega}|J-\tilde{J}|dl.
\end{equation}

Using this estimate and the coarea formula we have 
\begin{align*} 
\frac{m}{\sigma_1} \int_{\Omega} |u - \tilde{u}| dx & \leq  \int_{\R} \int_{\{u = t\}\cap \Omega} |u - \tilde{u}| dS dt \\
	& \leq  \frac{1}{c_g}\int_{\R} \int_{\R} \int_{\{u = t\} \cap \{g_t = r\} \cap \Omega}  |u - \tilde{u}| dl dr dt\\ 
	& \leq  \frac{L_M}{c_g \sigma_0}\int_{\R} \int_{\R} \int_{\{u = t\} \cap \{g_t = r\} \cap \Omega}  |J - \tilde{J}| dl dr dt\\
	 & =  \frac{L_M}{c_g \sigma_0} \int_{\R} \int_{\{u = t\}} |\nabla g_t | |J - \tilde{J}| dS dt  \\
	 & \leq \frac{L_M C_g}{c_g \sigma_0} \int_{\R} \int_{\{u = t\}}  |J - \tilde{J}| dS dt\\
	 \end{align*}	
\begin{align*} 
	 & =\frac{L_M C_g}{c_g \sigma_0} \int_{\Omega} |\nabla u| |J - \tilde{J}| dx  \\
	 & \leq \frac{L_M C_g}{c_g \sigma_0}  \|\nabla u \|_{L^{\infty} (\Omega)} \left(C \| |J| - |\tilde{J}| \|^{\frac{1}{2}}_{L^{\infty} (\Omega)} \right)\\
	 &\leq  \frac{L_M C_g M}{c_g \sigma_0^2} \left(C  \||J| - |\tilde{J}| \|^{\frac{1}{2}}_{L^{\infty} (\Omega)} \right),
\end{align*}
where we have applied  Theorem \ref{levelSetTheorem}. \hfill $\square$

\section{$W^{1,1}$ stability of the minimizers}

In this section, we prove stability of minimizers of \eqref{LGP} in $W^{1,1}$. As mentioned in Section 3, in general \eqref{LGP} does not even have unique minimizers, so in order to prove stability results in $W^{1,1}$,  it is natural to expect stronger assumptions on on the  minimizers. 

\begin{lemma}  Let $n=2,3$, and suppose $u$ and $\tilde{u}$ are admissible with $u|_{\partial \Omega}=\tilde{u}|_{\partial \Omega}=f\in L^{\infty}(\partial \Omega)$ and corrsponding conductivities $\sigma$ and $\tilde{\sigma}$, and current density vector fields $J$ and $\tilde{J}$, respectively. Suppose $\sigma, \tilde{\sigma} \in C^2(\bar{\Omega})$ with 
\begin{equation}\label{sigmaBound}
\parallel \sigma \parallel_{C^2(\Omega)}, \parallel \tilde{\sigma} \parallel_{C^2(\Omega)} \leq \sigma_2
\end{equation}
for some $\sigma_2 \in \R$. Let 
\begin{equation} \label{G}
G(x):=\frac{\tilde{J}(x) - J(x)}{\tilde{\sigma}(x)}, \ \ x\in \Omega,
\end{equation}
with $G=(G_1,G_2)$ for $n = 2$ and $G=(G_1,G_2, G_3)$ for $n = 3$. Then
\begin{equation}\label{GagliardoâNirenberg}
 \| \nabla G_i \|_{L^1 (\Omega)} \leq C_1 \|J-\tilde{J} \|^{1/2}_{L^1 (\Omega)},
\end{equation}
for some constant $C_1$ which depends only on $\Omega$, $\sigma_0$, $\sigma_2$ and $\parallel f \parallel_{L^{\infty}(\Omega)}$.
\end{lemma}
\vskip 1em \noindent
{\bf Proof.} Since $u$ and $\tilde{u}$ satisfy \eqref{PDE}, it follows from elliptic regularity that
\begin{equation} \label{ellipticRegularity}
\parallel u \parallel_{H^3(\Omega)}, \parallel \tilde{u} \parallel_{H^3(\Omega)}\leq C_1 \parallel f \parallel_{L^2(\Omega)}\leq C_1|\Omega| ^{\frac{1}{2}}\parallel f \parallel_{L^{\infty}(\Omega)},
\end{equation}
for some constant $C_1$ depending only on $\sigma_0$, $\sigma_2$, and $\Omega$. Now note that 
\[G(x)=\nabla \tilde{u}-\frac{\sigma}{\tilde{\sigma}}\nabla u.\]
Thus it follows from  \eqref{sigmaBound} and \eqref{ellipticRegularity} that 
\begin{equation}\label{L1L2}
\parallel D^2 G_i \parallel_{L^1(\Omega)} \leq |\Omega|^{\frac{1}{2}}\parallel D^2 G_i \parallel_{L^2(\Omega)} \leq C,\ \  1 \leq i\leq n, 
\end{equation}
for some constant $C$ which only depends on $\sigma_0$, $\sigma_2$, $\Omega$ and $\parallel f \parallel_{L^{\infty}(\Omega)}$. On the other hand it follows from Gagliardo-Nirenberg interpolation inequality that
\begin{equation}\label{Gagliardo-Nirenberg}
 \| \nabla G_i \|_{L^1 (\Omega)} \leq C_2 \| D^2G_i \|^{1/2}_{L^1 (\Omega)} \|G_i \|^{1/2}_{L^1 (\Omega)},
\end{equation}
for some $C_2$ which only depends on $\Omega$. Combining \eqref{L1L2},  \eqref{Gagliardo-Nirenberg}, and 
\[\parallel G_i \parallel_{L^1(\Omega)}\leq \frac{\parallel J-\tilde{J} \parallel_{L^1(\Omega)}}{\sigma_0}, 1 \leq i\leq n,\]
we arrive at the inequality \eqref{Gagliardo-Nirenberg}. \hfill $\square$

\vskip 1em \noindent
Next we prove that $u$ and $\tilde{u}$ are close in $W^{1,1}(\Omega)$. In order to do so, we need additional assumptions on the structure of level sets of $u$. 

\begin{definition} \label{WellStructured} Suppose $u$ is admissible, $n=2$, and $x\in \Omega$. Pick $h \in \R^2$ with $|h| = 1$, and $t \in \R$ small enough such that $x + th \in \Omega$. Let $\Gamma$ and $\Gamma_t$ be the level sets of $u$ passing through $x$ and $x+th$, respectively. Parametrize $\Gamma$ and $\Gamma_t$ by the arc length such that $\gamma(0), \gamma_t(0)\in \partial \Omega$, and denote these parametrizations by $\gamma$ and $\gamma_t$, respectively. 

Similarly in dimension $n=3$, let  $u$ be admissible and suppose level sets of $u$ can be foliated to one-dimensional curves in the sense of Definition \ref{foliation}. 
Pick $x \in \Omega$ and $h\in \R^3$ with $|h|=1$, and choose $t$ small enough such that $x+th \in \Omega$. Let $\Gamma$ and $\Gamma_t$ be the unique curves in 
\[\{\{u=\tau\} \cap \{g_{\tau}=r\} \ \  \tau, r \in \R\}\]
which pass through $x$ and $x+th$, respectively, and let  $\gamma$ and $\gamma_t$ be the parametrization of these curves with respect to arc length with $\gamma(0), \gamma_t(0)\in \partial \Omega$. 
 
 We say that level sets of $u$ are \emph{well structured}  if the following conditions are satisfied 
\begin{enumerate}[label=(\alph*)]

\item There exists $K \geq 0$ such that
\begin{equation}\label{(a)}
 \left| \frac{\gamma_{t}^{'} (s) - \gamma^{'} (s)}{t} \right| \leq K
\end{equation}
for every $s \in [0,L]$, $t\in \R$, $x\in \Omega$ and $h\in S^
{n-1}$. In particular,
\begin{equation} \label{(b)}
\gamma_{t}^{'} (s) \to \gamma^{'} (s) \ \ \hbox{as}\ \  t \to 0,
\end{equation}
where $\gamma'(s)=\frac{d \gamma(s)}{ds}$ and $\gamma'_t(s)=\frac{d \gamma_t(s)}{ds}$. 
\item There exists a bounded function $F_{x,h} (s)=F(x,h;s) \in L^{\infty} (\Omega \times S^{n-1} \times [0,L_M])$ such that
\begin{equation}\label{(c)}
 \lim_{t \to 0} \frac{\gamma_{t} (s) - \gamma (s)}{t} = F_{x,h}(s)
\end{equation}
for every $s \in [0,L]$, $x\in \Omega$ and $h\in S^{n-1}$.
\end{enumerate}

\end{definition}

\begin{remark}
Let $x\in \Omega$, $h \in \R^2$ with $|h| = 1$, and $t \in \R$ be small enough such that $x + th \in \Omega$. Also, as in Definition \ref{WellStructured}, let $\gamma$, and $\gamma_t$ be the parametrization of the curves passing through $x$ and $x+th$. In view of Remark \ref{foliationRemark} we have 

\begin{equation}
\gamma(s)=\mathcal{F}_{u(x)}(\bar{\gamma}(s)) \ \ \ \ \hbox{and}\ \ \ \  \gamma_t(s)=\mathcal{F}_{u(x+th)}(\bar{\gamma}_t(s)),
\end{equation}
where $\bar{\gamma}(s)$ and $\bar{\gamma}_t(s)$ are parametrization of two level sets of the function $g(x,y)=y=\Pi_y(\mathcal{F}^{-1}(x))$ and  $g(x,y)=y=\Pi_y(\mathcal{F}^{-1}(x+th))$, respectively. Here $\Pi_y$ is the projection operator on $y$-axis, and $\mathcal{F}_{u(x)}$ and $\mathcal{F}_{u(x+th)}$ are $C^1$ diffeomorphisms from $B(0,1)$ to the connected components of  the level sets of $u$ passing through $x$ and $x+th$, respectively.  It is easy to see that $\bar{\gamma}_t(s)$ is continuously differentiable with respect to $t$, for each fixed $s$.

Now let $\Gamma_{x_0}$ be the connected component of the level set of $u$ that passes through $x_0$, and assume that $|\nabla u|>0$ on $\Omega$. Then in a neighborhood of $r_0=u(x_0)$ we can find $C^1$ diffeomorphisms $F_r$ so that $F_r (y)$ is continuously differentiable with respect to $r$, for each fixed $y$. Indeed let $y\in B(0,1)$ and consider the gradient flow 
\begin{equation}
\dot{z}_y(q)=\nabla u(z_y(q)), \ \ z_y(0)=F_0(y),
\end{equation} 
which has a unique solution as long as $z_y(q) \in \Omega$. Let $r \in range(u)$ be and $\Gamma_r$ be a connected component of $\{u=r\}$. Define $F_r: B(0,1) \rightarrow \Gamma_r$ by 
\[F_r(y)=F_{r_0}(z_y(q_r)),\]
where $q_r \in \mathbb{R}$ is the unique point where $z_y(q_r)\in \Gamma_r$. Also observe that the set 
\[\mathcal{R}=\{ \ r\in range(u): \ \  \mathcal{F} \ \ \hbox{is well defined on} \ \ \{u=r\} \},\]
is both open and closed in $range(u)$, and hence $\mathcal{R}=range(u)$ and therefore $F_r$ could be defined globally as above for all $r\in range(u)$.  

Since $u$, $F_{r_0}$, and $z_y$ are all $C^1$, it is easy to see that $F_r(y)$ is continuously differentiable with respect to $r$, for each fixed $y\in B(0,1)$. Now notice that the level sets of the function $g(x,y): B(0,1) \rightarrow \mathbb{R}$ defined by $g(x,y)=y$ are well structured in the sense of Definition \ref{WellStructured}. In view of the above arguments, it follows from the chain rule that $\gamma_t(s)=\mathbb{F}_t(\bar{\gamma}_t(s))$, where $\bar{\gamma}_t(s)$ is a parametrization of the level set $g(x,y)=y$ passing through $\mathcal{F}^{-1}(x+th)$, and $\mathcal{F}_t$ and $\bar{\gamma}_t$ are both continuously differentiable with respect to $t$. Therefore, since $\eqref{(a)}, \eqref{(b)}, \eqref{(c)}$ hold for any parametrization of level sets of $g(x,y)=y$, an application of the chain rule implies that $\eqref{(a)}, \eqref{(b)}, \eqref{(c)}$ also hold under the assumptions of Definition \eqref{WellStructured}. In particular, if $u$ is a $C^1$ function with $|\nabla u|>0$ in $\Omega$ and $\{x\in \partial \Omega: f(x)=t\}$ has finitely many connected components for all $t$, then level sets of $u$ are well structured in the sense of Definition \ref{WellStructured}. 

\end{remark}

\vskip 1em \noindent
\begin{theorem}  \label{gradientcloseTheo} Let $n=2$, and suppose $u$ and $\tilde{u}$ are admissible with $u|_{\partial \Omega}=\tilde{u}|_{\partial \Omega}=f,$  corresponding conductivities $\sigma, \tilde{\sigma} \in C^2(\Omega)$, and current density vector fields $J$ and $\tilde{J}$, respectively. Suppose $\sigma, \tilde{\sigma} \in C^2(\bar{\Omega})$ and satisfy \eqref{sigmaBound}. If $u$ satisfies \eqref{bpundedLength}, and the level sets of $u$ are well-structured in the sense of Definition \ref{WellStructured}, then 

\begin{equation}\label{gradientsAreClose}
 \| \nabla \tilde{u} - \nabla u \|_{L^1 (\Omega)} \leq C\parallel |J| - |\tilde{J}|\parallel^{\frac{1}{4}}_{L^{\infty}({\Omega})},
\end{equation}
for some constant $C(m,M,\sigma_0,\sigma_1,\sigma_2, u, f, L_M)$ independent of $\tilde{u}$ and $\tilde{\sigma}$. 
\end{theorem}

\vskip 1em \noindent
\textbf{Proof.} Fix $x\in \Omega$ and $h \in \R^2$ with $|h|=1$. Then 
\[  \mathcal{L}(x,h):=(\nabla \tilde{u}(x) - \nabla u(x)) \cdot h = \lim_{t \to 0} \frac{[\tilde{u}(x + t h) - u(x + t h)] - [\tilde{u}(x) - u(x)]}{t}. \]
First we estimate the above limit. Since all level sets of $u$ reach the boundary $\partial \Omega$, there exist $z,z_t \in \partial \Omega$ such that
\[  u(x) = u(z) = \tilde{u} (z),  \]
\[  u(x + t h) = u(z_t) = \tilde{u} (z_t). \]
Thus
\[  [\tilde{u}(x + t h) - u(x + t h)] - [\tilde{u}(x) - u(x)] = [\tilde{u}(x + t h) - \tilde{u}(z_t)] - [\tilde{u}(x) - \tilde{u}(z)]. \]

Let $\gamma$ and $\gamma_t$ be the curves passing through $x$ and $x+th$, described in Definition \ref{WellStructured} with $\gamma(0)=z$ and $\gamma_t(0)=z_t$. Suppose $\gamma(s_0)=x$ and reparamterize $\gamma_t$ so that $\gamma_t(s_0)=x+th$. Then we have
\[   [\tilde{u}(x + t h) - \tilde{u}(z)] - [\tilde{u}(x) - \tilde{u}(z)] = [\tilde{u}(\gamma_t (s_0)) - \tilde{u}(\gamma_t (0))] - [\tilde{u}(\gamma(s_0)) - \tilde{u}(\gamma(0))]  \]
\[   = \int_{0}^{s_0} \nabla \tilde{u} (\gamma_t (s)) \cdot \gamma'_t (s) ds - \int_{0}^{s_0} \nabla \tilde{u} (\gamma (s)) \cdot \gamma' (s) ds.   \]
Hence 
\[   \mathcal{L}(x,h)= \lim_{t \to 0} \frac{1}{t} \left( \int_{0}^{s_0} \nabla \tilde{u} (\gamma_t (s)) \cdot \gamma'_t (s) ds - \int_{0}^{s_0} \nabla \tilde{u} (\gamma (s)) \cdot \gamma' (s) ds \right). \]
Substituting $\nabla \tilde{u}$ by $\frac{\tilde{J}}{\tilde{\sigma}}$ and using the fact that $J$ is perpendicular to $\gamma'$ and $\gamma'_t$ we get
\[   \mathcal{L}(x,h)=\lim_{t \to 0} \frac{1}{t} \left( \int_{0}^{s_0} \frac{\tilde{J}(\gamma_t (s)) - J(\gamma_t (s))}{\tilde{\sigma}(\gamma_t (s))} \cdot \gamma'_t (s) ds - \int_{0}^{s_0} \frac{\tilde{J}(\gamma (s)) - J(\gamma(s))}{\tilde{\sigma}(\gamma (s))} \cdot \gamma' (s) ds \right). \]
Now define 
\[  G(x):=\frac{\tilde{J}(x) - J(x)}{\tilde{\sigma}(x)}, \ \ x\in \Omega. \]
Hence 
\[   \mathcal{L}(x,h)=\lim_{t \to 0} \frac{1}{t} \left ( \int_{0}^{s_0} G(\gamma_t (s)) \cdot \gamma'_t (s) ds - \int_{0}^{s_0} G(\gamma (s)) \cdot \gamma' (s) ds \right). \]
The expression in the right hand side can be rewritten as
\begin{equation}\label{limit}
\frac{1}{t}\int_{0}^{s_0} [G(\gamma_t (s)) - G(\gamma (s)) ] \cdot \gamma'_t (s) ds + \frac{1}{t}\int_{0}^{s_0} G(\gamma (s)) \cdot [\gamma'_t (s) - \gamma' (s) ]ds.
\end{equation}
It follows from the assumption (a) in Definition \ref{WellStructured} that 
\[  \left| \frac{\gamma'_t (s) - \gamma' (s)}{t} \right| \leq K, \]
and hence 
\begin{equation}\label{secondTerm}
\left| \frac{1}{t} \int_{0}^{s_0} G(\gamma (s)) \cdot [\gamma'_t (s) - \gamma' (s) ]ds \right| \leq \frac{K}{\sigma_0} \int_{0}^{L}  |\tilde{J}(\gamma (s)) - J(\gamma(s))| ds.
\end{equation}

\vskip 1em \noindent
Now we turn our attention to the first term in \eqref{limit}. Let $G=(G_1,G_2)$. Since
\[ \lim_{t \to 0} \frac{\gamma_t (s) - \gamma (s)}{t}=F_{x,h}(s)\]
for $i=1,2$ we have 
\begin{eqnarray*}
 \lim_{t \to 0} \frac{G_i(\gamma_t (s)) - G_i(\gamma (s))}{t}&=&\lim_{t \to 0} \frac{G_i(\gamma (s)+tF(s)) - G_i(\gamma (s))}{t} \\
 &=&  \nabla G_i (\gamma (s)) \cdot F(s). \\ 
\end{eqnarray*}
Thus the first term of \eqref{limit} can be rewritten as
\begin{eqnarray}\label{firstterm}
    && \lim_{t \to 0} \frac{1}{t}  \int_{0}^{s_0} [G(\gamma_t (s)) - G(\gamma (s)) ] \cdot \gamma'_t (s) dl \nonumber \\
    && = \int_{0}^{s_0} (\nabla G_1 (\gamma (s)) \cdot F(s), \nabla G_2 (\gamma (s)) \cdot F(s)) \cdot \gamma' (s) dl \nonumber \\
    && \leq \parallel F \parallel_{L^\infty} \int_{0}^{s_0} | \nabla G_1 (\gamma (s)) | +| \nabla G_2 (\gamma (s)) | dl \nonumber \\
    && \leq \parallel F \parallel_{L^\infty} \int_{0}^{L} | \nabla G_1 (\gamma (s)) | +| \nabla G_2 (\gamma (s)) | dl,
\end{eqnarray}
where we have used the assumption (b) in Definition \ref{WellStructured}. Combining \eqref{secondTerm} and \eqref{firstterm} we obtain

\begin{eqnarray*}
| \nabla \tilde{u}(x)-\nabla u(x)| & \leq & \sup_{h\in \R^2,  |h|=1} \mathcal{L}(x,h)  \\
&\leq & \frac{K}{\sigma_0} \int_{0}^{L}  |\tilde{J}(\gamma (s)) - J(\gamma(s))| dl\\
&&+ \parallel F \parallel_{L^\infty} \int_{0}^{L} | \nabla G_1 (\gamma (s)) | +| \nabla G_2 (\gamma (s)) | dl.
\end{eqnarray*}
Thus
\begin{eqnarray*}
\int_{\Gamma} | \nabla \tilde{u}(x)-\nabla u(x)| dl &\leq & \frac{KL_M}{\sigma_0} \int_{\Gamma}  |\tilde{J}(x) - J(x)| dl\\
&&+ L_M\parallel F \parallel_{L^\infty} \int_{\Gamma} | \nabla G_1 (x) | +| \nabla G_2 (x) | dl,\nonumber
\end{eqnarray*}
and consequently 
\begin{eqnarray}\label{firstAndSecond}
\int_{\{u=\tau\}\cap \Omega} | \nabla \tilde{u}(x)-\nabla u(x)| dl &\leq & \frac{KL_M}{\sigma_0} \int_{\{u=\tau\}\cap \Omega}  |\tilde{J}(x) - J(x)| dl\\
&&+ L_M\parallel F \parallel_{L^\infty} \int_{\{u=\tau\}\cap \Omega} | \nabla G_1 (x) | +| \nabla G_2 (x) | dl.\nonumber
\end{eqnarray}

\vskip 1em \noindent
Using \eqref{firstAndSecond} and the coarea formula we have 
\begin{eqnarray*}
 \frac{m}{\sigma_1} \| \nabla \tilde{u} - \nabla u \|_{L^1 (\Omega)} & \leq & \int_{\Omega} |\nabla u| |\nabla \tilde{u} - \nabla u| dx \\
     & =  &\int_{\R} \int_{\{u = \tau\}\cap \Omega} |\nabla \tilde{u} - \nabla u| dl d\tau \\
     & \leq & \frac{KL_M}{\sigma_0} \int_{\R} \int_{\{u=\tau\}\cap \Omega}  |\tilde{J} - J| dl d\tau\\
&+& L_M\parallel F \parallel_{L^\infty} \int_{\R} \int_{\{u=\tau\}\cap \Omega} | \nabla G_1 | +| \nabla G_2 | dl d\tau \\
& \leq & \frac{KL_M M}{ (\sigma_0)^2} \int_{\R} \int_{\{u=\tau\}\cap \Omega} \frac{ |\tilde{J}- J| }{ |\nabla u|}dl d\tau\\
&+&\frac{ L_M\parallel F \parallel_{L^\infty }M}{\sigma_0} \int_{\R} \int_{\{u=\tau\}\cap \Omega}  \frac{ |\nabla G_1  | +| \nabla G_2  | }{ |\nabla u|} dl d\tau \\
& = & \frac{KL_M M}{ (\sigma_0)^2} \int_{ \Omega}  |\tilde{J} - J| dx\\
&+&\frac{ L_M\parallel F \parallel_{L^\infty }M}{\sigma_0} \int_{\Omega}  |\nabla G_1 | +| \nabla G_2  |  dx\\
&\leq & \frac{KL_M M}{ (\sigma_0)^2} \parallel J - \tilde{J}  \parallel_{L^1(\Omega)}\\
&+&\frac{2 L_M C_1 \parallel F \parallel_{L^\infty }M}{\sigma_0} \parallel J-\tilde{J}\parallel^{\frac{1}{2}}_{L^1(\Omega)}\\
\end{eqnarray*}
where we have used \eqref{Gagliardo-Nirenberg} to obtain the last inequality. Applying Theorem \ref{levelSetTheorem}, and noting that 
\[ \parallel J-\tilde{J}\parallel^{\frac{1}{2}}_{L^1(\Omega)} \leq (2M|\Omega|)^{\frac{1}{2}},\]
where $M$ is defined in \eqref{mM}, we arrive at  \eqref{gradientsAreClose}.  
\hfill $\square$
\vskip 1em \noindent
Now we prove three dimensional version of this theorem.
\begin{theorem} \label{gradientAreCloseTheo1}
Let $n=3$, and suppose $u$ and $\tilde{u}$ are admissible with $u|_{\partial \Omega}=\tilde{u}|_{\partial \Omega}=f,$  corresponding conductivities $\sigma, \tilde{\sigma} \in C^2(\Omega)$, and current density vector fields $J$ and $\tilde{J}$, respectively. Suppose $\sigma, \tilde{\sigma} \in C^2(\bar{\Omega})$ and satisfy \eqref{sigmaBound}. In addition suppose $u$ satisfies \eqref{bpundedLength}, the level sets of $u$ can be foliated to one-dimensional curves in the sense of Definition 3.4, and the level sets of $u$ are well-structured in the sense of Definition 4.2. Then 

\begin{equation}\label{gradientsAreClose1}
 \| \nabla \tilde{u} - \nabla u \|_{L^1 (\Omega)} \leq C\|a - \tilde{a}\|^{\frac{1}{4}}_{L^{\infty}({\Omega})},
\end{equation}
for some constant $C(m,M,\sigma_0,\sigma_1,\sigma_2, u, f, L_M, c_g,C_g)$ is independent of $\tilde{u}$ and $\tilde{\sigma}$. 

\end{theorem}

\vskip 1em \noindent
\textbf{Proof.} With an argument similar to the one used in the proof of Theorem \ref{gradientcloseTheo} we get 
\begin{eqnarray}\label{firstAndSecond1}
\int_{U_{\tau,r}} | \nabla \tilde{u}(x)-\nabla u(x)| dl &\leq & \frac{KL_M}{\sigma_0} \int_{U_{\tau,r}}  |\tilde{J}(x) - J(x)| dl\\
&&+ L_M\parallel F \parallel_{L^\infty} \int_{U_{\tau,r}} | \nabla G_1 (x) | +| \nabla G_1 (x) | + | \nabla G_3 (x) | dl,\nonumber
\end{eqnarray}
where  $U_{\tau,r} := \{u=\tau\}\cap\{g_{\tau}=r\}\cap \Omega$  and $G=(G_1,G_2,G_3)$ is defined in \eqref{G}.

It follows follows from \eqref{firstAndSecond1} and  the coarea formula that 
\begin{eqnarray*}
 \frac{m}{\sigma_1} \| \nabla \tilde{u} - \nabla u \|_{L^1 (\Omega)} & \leq & \int_{\Omega} |\nabla u| |\nabla \tilde{u} - \nabla u| dx \\
     & =  &\int_{\R} \int_{\{u = \tau\}\cap \Omega} |\nabla \tilde{u} - \nabla u| dS d\tau \\
     & =  &\int_{\R} \int_{\{u = \tau\}\cap \Omega} \frac{|\nabla g_{\tau}|}{|\nabla g_{\tau}|} |\nabla \tilde{u} - \nabla u| dS d\tau \\
     & =  &\int_{\R} \int_{\R} \int_{U_{\tau,r}} \frac{1}{|\nabla g_{\tau}|} |\nabla \tilde{u} - \nabla u| dl dr d\tau \\
     & \leq & \frac{KL_M}{\sigma_0 c_g} \int_{\R} \int_{\R} \int_{U_{\tau,r}}  |\tilde{J} - J| dl dr dt\\ 
&+& \frac{L_M\parallel F \parallel_{L^\infty}}{c_g} \int_{\R} \int_{\R} \int_{U_{\tau,r}} | \nabla G_1 | +| \nabla G_2 | + | \nabla G_3 | dl dr dt \\   
& \leq & \frac{KL_M M C_g}{ (\sigma_0)^2 c_g} \int_{\R} \int_{\R} \int_{U_{\tau,r}} \frac{ |\tilde{J}- J| }{ |\nabla u| |\nabla g_{\tau}|}dl dr dt\\
&+&\frac{ L_M M\parallel F \parallel_{L^\infty } C_g}{\sigma_0 c_g} \int_{\R} \int_{\R} \int_{U_{\tau,r}}  \frac{ |\nabla G_1  | +| \nabla G_2  | + | \nabla G_3  |}{ |\nabla u| |\nabla g_t|} dl dr dt \\
& = & \frac{KL_M M C_g}{ (\sigma_0)^2 c_g} \int_{\R} \int_{\{u = \tau\}\cap \Omega} \frac{ |\tilde{J}- J| }{ |\nabla u| }dS dt\\
&+&\frac{ L_M M\parallel F \parallel_{L^\infty} C_g}{\sigma_0 c_g} \int_{\R} \int_{\{u = \tau\}\cap \Omega}  \frac{ |\nabla G_1  | +| \nabla G_2  | + | \nabla G_3  |}{ |\nabla u|} dS dt \\
& = & \frac{KL_M M C_g}{ (\sigma_0)^2 c_g} \int_{ \Omega}  |\tilde{J} - J| dx\\
&+&\frac{ L_M M\parallel F \parallel_{L^\infty} C_g}{\sigma_0 c_g} \int_{\Omega}  |\nabla G_1 | +| \nabla G_2  | + | \nabla G_3  | dx\\
 \end{eqnarray*}    
 \begin{eqnarray*}
&\leq & \frac{KL_M M C_g}{ (\sigma_0)^2} \parallel J - \tilde{J}  \parallel_{L^1(\Omega)}\\
&+&\frac{2 L_M C_1 M \parallel F \parallel_{L^\infty(\Omega) }C_g}{\sigma_0} \parallel J-\tilde{J}\parallel^{\frac{1}{2}}_{L^1(\Omega)},
\end{eqnarray*}
where we have used \eqref{Gagliardo-Nirenberg} to obtain the last inequality. Applying Theorem \ref{levelSetTheorem}, and noting that 
\[ \parallel J-\tilde{J}\parallel^{\frac{1}{2}}_{L^1(\Omega)} \leq(2M|\Omega|)^{\frac{1}{2}},\]
we obtain the inequality  \eqref{gradientsAreClose}.  
\hfill $\square$

\vskip 1em \noindent
Now, we are ready to prove our main stability results. 
\begin{theorem} \label{sigmaStability2}
Let $n=2$, and suppose $u$ and $\tilde{u}$ are admissible with $u|_{\partial \Omega}=\tilde{u}|_{\partial \Omega}=f,$  corresponding conductivities $\sigma, \tilde{\sigma} \in C^2(\Omega)$, and current density vector fields $J$ and $\tilde{J}$, respectively. Suppose $\sigma, \tilde{\sigma} \in C^2(\bar{\Omega})$ and satisfy \eqref{sigmaBound}. If $u$ satisfies \eqref{bpundedLength} and level sets of $u$ are well-structured in the sense of Definition \ref{WellStructured},  then
\begin{equation*}
    \|  \sigma - \tilde{\sigma} \|_{L^1 (\Omega)} \leq C\parallel |J| - |\tilde{J}| \parallel^{\frac{1}{4}}_{L^{\infty}({\Omega})},
\end{equation*}
for some constant $C(m,M,\sigma_0,\sigma_1,\sigma_2, \sigma, f, L_M)$  independent of $\tilde{\sigma}$. 
\end{theorem}

\vskip 1em \noindent
\textbf{Proof.} Using Theorem \ref{gradientcloseTheo} we have 
\begin{eqnarray*}
    \int_{\Omega} | \sigma - \tilde{\sigma} | dx & = & \int_{\Omega} \left| \frac{|J|(|\nabla \tilde{u}| - |\nabla u|)}{|\nabla u||\nabla \tilde{u}|} + \frac{|J| - |\tilde{J}|}{|\nabla \tilde{u}|} \right| dx \\
& \leq & \int_{\Omega}  \frac{|J|}{|\nabla u||\nabla \tilde{u}|} \left| |\nabla u| - |\nabla \tilde{u}| \right| dx 
+ \int_{\Omega} \frac{1}{|\nabla \tilde{u}|} \left| |J| - |\tilde{J}| \right| dx \\
& \leq & \int_{\Omega}  \frac{|J|}{|\nabla u||\nabla \tilde{u}|}  |\nabla u - \nabla \tilde{u}|  dx 
+ \int_{\Omega} \frac{1}{|\nabla \tilde{u}|} \left| |J| - |\tilde{J}| \right| dx \\
& \leq & \frac{M \sigma_{1}^{2} C}{m^2} \parallel |J| - |\tilde{J}| \parallel^{\frac{1}{4}}_{L^{\infty}(\Omega)} 
+ \frac{\sigma_1 |\Omega|}{m} \parallel |J| - |\tilde{J}| \parallel_{L^{\infty}(\Omega)} \\
& \leq & \left[ \frac{M \sigma_{1}^{2} C}{m^2} + \frac{\sigma_1 |\Omega| (2M)^{\frac{3}{4}}}{m} \right] \parallel |J| - |\tilde{J}| \parallel^{\frac{1}{4}}_{L^{\infty}(\Omega)}. \\
\end{eqnarray*}
\hfill $\Box$

\begin{theorem}\label{sigmaStability3}
Let $n=3$, and suppose $u$ and $\tilde{u}$ are admissible with $u|_{\partial \Omega}=\tilde{u}|_{\partial \Omega}=f,$  corresponding conductivities $\sigma, \tilde{\sigma} \in C^2(\Omega)$, and current density vector fields $J$ and $\tilde{J}$, respectively. Suppose $\sigma, \tilde{\sigma} \in C^2(\bar{\Omega})$ and satisfy \eqref{sigmaBound}. If $u$ satisfies \eqref{bpundedLength}, the level sets of $u$ can be foliated to one-dimensional curves in the sense of Definition \ref{foliation}, and the level sets of $u$ are well-structured in the sense of Definition \ref{WellStructured}, then 

\begin{equation}
 \| \sigma - \tilde{\sigma}\|_{L^1 (\Omega)} \leq C\||J| - |\tilde{J}|\|^{\frac{1}{4}}_{L^{\infty}({\Omega})},
\end{equation}
for some constant $C(m,M,\sigma_0,\sigma_1,\sigma_2, \sigma, f, L_M, c_g,C_g)$ independent of $\tilde{\sigma}$.

\end{theorem}
\vskip 1em \noindent
\textbf{Proof.} The proof follows from Theorem \ref{gradientAreCloseTheo1} and a calculation similar to that of the proof of Theorem \ref{sigmaStability2}. $\square$


\begin{thebibliography}{99}

\bibitem{Alberti}
G. Alberti
{\it A Lusin type theorem for gradients}, 
J. Funct. Anal. 100 (1991), no. 1, 110-118.

\bibitem{Ales} G. Alessandrini, 
\textit{Critical points of solutions of elliptic equations in two variables.}
Ann. Scuola Norm. Sup. Pisa Cl. Sci. (4) 14 (1987), no. 2, 229-256 (1988).


\bibitem{Anzellotti}
G. Anzellotti, 
{\it Pairings between measures and bounded functions and compensated compactness}, 
Ann. Mat. Pura Appl. (4) 135 (1983), 293-318 (1984). 

\bibitem{borcea}{ L. Borcea}, {\em Electrical impedance tomography}, Inverse Problems {\bf 18}(2002), R99-R136.
\bibitem{CI} M. Cheney and D. Isaacso, \textit{An overview of inversion algorithms for impedance imaging}, Contemp. Math. (1991) 122, 29-39.
\bibitem{isaacsonReview}{ M. Cheney, D. Isaacson, and J. C. Newell}, {\em Electrical Impedance Tomography},
SIAM Rev. {\bf 41}(1999), no.1, 85-101.

\bibitem{Dos} M. Dos Santos, \textit{Characteristic functions on the boundary of a planar domain need not be traces of least gradient functions},  Confluentes Math. \textbf{9} (2017), no. 1, 65-93.

\bibitem{G} W. G\'{o}rny, P. Rybka, and A. Sabra, \textit{Special cases of the planar least gradient problem}.
Nonlinear Anal. \textbf{151} (2017), 66-95.

\bibitem{GLKU}{ A. Greenleaf, Y. Kurylev, M. Lassas, and G. Uhlmann}, {\em Invisibility and
inverse problems}, Bull. Amer. Math. Soc. 46 (2009), 55-97.

\bibitem {G} E. Giusti, \textit{Minimal Surfaces and Functions of Bounded Variations}, 1984 (Boston:
Birkh\"a user).

\bibitem{HMN} N. Hoell, A. Moradifam, A. Nachman, \textit{Current density impedance imaging of an anisotropic conductivity in a known conformal class}, SIAM J. Math. Anal. 46 (2014), no. 3, 1820-1842.

\bibitem{isacson} {\sc D. Isaacson and M. Cheney},
{\em Effects of measurement precision and finite numbers of
electrodes on linear impedance imaging algorithms,} SIAM J. Appl.
Math. 51 (1991), no. 6, 1705-1731.
\bibitem{JMN} R.L. Jerrard, A. Moradifam, A.I. Nachman. \textit{Existence and uniqueness of minimizers of general least gradient problems.} J. Reine Angew. Math. 734 (2018), 71-97.

\bibitem{joy89} { M. L. Joy, G. C. Scott, and M.
Henkelman}, {\em In vivo detection of applied electric currents by magnetic resonance imaging}, Magnetic Resonance Imaging, {\bf 7} (1989), pp. 89-94.
\bibitem{joy}
{\sc M. J. Joy, A. I. Nachman, K. F. Hasanov, R. S. Yoon, and A.
W. Ma}, {\em A new approach to Current Density Impedance Imaging
(CDII)}, Proceedings ISMRM, No. 356, Kyoto, Japan, 2004.
\bibitem{Man} N. Mandache, Exponential instability in an inverse problem for the Schrodinger equation Inverse Problems (2001) \textbf{17} 1435-1444.
\bibitem{MT} C. Montalto, A. Tamasan, \textit{Stability in conductivity imaging from partial measurements of one interior current}, Inverse Probl. Imaging \textbf{11} (2017), no. 2, 339-353. 

\bibitem{MS} C. Montalto, P. Stefanov,  \textit{Stability of coupled-physics inverse problems with one internal measurement}, Inverse Problems \textbf{29} (2013), no. 12, 125004, 13 pp. 

\bibitem{Struc2} A. Moradifam,  \textit{Least gradient problems with Neumann boundary condition}, Journal of Differential Equations, 263 (2017), no. 11, 7900-7918.
\bibitem {Struc1} A. Moradifam, \textit{Existence and structure of minimizers of least gradient problems}, Indiana University Math Journal, 67 No. 3 (2018), 1025-1037. 
\bibitem{MNT} A. Moradifam,  A. Nachman, A. Tamasan,  Conductivity imaging from one interior measurement in the presence of perfectly conducting and insulating inclusions, SIAM J. Math. Anal., 44 (2012), 3969-3990

\bibitem{MNTim} A. Moradifam, A. Nachman, A. Timonov, \textit{A convergent algorithm for the hybrid problem of reconstructing conductivity from minimal interior data}, Inverse Problems, \textbf{28} (2012) 084003 (23pp). 

\bibitem{MNTCalc} A. Moradifam, A. Nachman, A. Tamasan, \textit{Uniqueness of minimizers of weighted least gradient problems arising in hybrid inverse problems}, Calc. Var. Partial Differential Equations \textbf{57} (2018), no. 1, Art. 6, 14 pp. 

\bibitem{Morris}  M. W. Hirsch, \textit{Differential topology}. Graduate Texts in Mathematics, No. 33. Springer-Verlag, New York-Heidelberg, 1976.

\bibitem{NTT07}{ A. Nachman, A. Tamasan, and A. Timonov}, {\em Conductivity
imaging with a single measurement of boundary and interior data},
Inverse Problems, {\bf 23} (2007), pp. 2551--2563.

\bibitem{NTT08}
{ A. Nachman, A. Tamasan, and A. Timonov}, {\em Recovering the
conductivity from a single measurement of interior data}, Inverse
Problems, {\bf 25} (2009) 035014 (16pp).

\bibitem{NTT10}{ A. Nachman, A. Tamasan, and A. Timonov}, {\em Reconstruction of Planar Conductivities in Subdomains from Incomplete Data},
SIAM J. Appl. Math. {\bf 70}(2010), Issue 8, pp. 3342--3362.

\bibitem{NTT11}{ A. Nachman, A. Tamasan, and A. Timonov}, {\em Current density impedance imaging}, Tomography and inverse transport theory, 135-149, Contemp. Math., 559, Amer. Math. Soc., Providence, RI, 2011.
\bibitem{nashedTa10}
{ M. Z. Nashed and A. Tamasan}, {\em Structural stability in a
minimization problem and applications to conductivity imaging},
Inverse Probl. Imaging, {\bf 4} (2010) to appear.

\bibitem{ST}
{ G.S. Spradlin and A. Tamasan}, {\em Not all traces on the circle come from functions of least gradient in the disk},
Indiana Univ. Math. J. \textbf{63} (2014), no. 6, 1819-1837.


\bibitem{Zuniga} A. Zuniga, \textit{Continuity of minimizers to weighted least gradient problems.} 
\textit{Nonlinear Anal.} 178 (2019), 86-109.



\end{thebibliography}
\end{document}